**Note On The Irrationality of Special Values of the *L*-Function**
**N. A. Carella, October, 2012.**

*Abstract:* A unified proof of the irrationality of the special values $L(n,\chi)$, $n \geq 2$, of beta *L*-function is put forward in this note. The first case of $n = 2$ seems to confirm that the Catalan constant $L(2,\chi)$ is an irrational number.

# 1 Introduction

Although most real numbers are irrational numbers, the classification of numbers as either rational or irrational numbers is an intractable problem. As the various known irrationality criteria are difficult to apply, the rationality or irrationality nature of unknown numbers are usually handled on a case-by-case basis. Many prominent constants such as the first irrational $\sqrt{2}$, the logarithmus naturalis base $e = \lim_{x \to n}(1 + n^{-1})^n = 2.7182818...$, the circle constant $\pi = 3.141592...$, the logarithm constant log 2, the special zeta values $\zeta(2)$, $\zeta(3)$, ..., et cetera, have their own tailor made proof of irrationality. An expanding collection of irrational numbers has been proven to be irrationals, and sometimes transcendental numbers using various techniques. In contrast, the rationality or irrationality nature of many other ubiquitous constants remains unknown, see [FH] for an extensive list of constants in the mathematical sciences. Some of the constants whose arithmetic properties are of current research interest are

1) Artin constant $C_A = \prod_{p \geq 3}(1 - p^{-1}(p-1)^{-1}) = 0.3739558...$,  (1)
2) Catalan constant $L(2,\chi) = \sum_{n \geq 1}\chi(n)n^{-2} = 0.915\,965\,594\,...$,
3) Euler constant $\gamma = \lim_{x \to \infty}\sum_{n \leq x} n^{-1} - \log x = 0.57721566...$,
4) Khintchin constant $C_K = \lim_{n \to \infty}\prod_{n \geq 0}(a_0 a_1 a_2 \cdots a_n)^{1/n} = 2.68545200...$, the geometric average of the continued fraction $[a_0, a_1, ...]$ of a real number,
5) Landau-Ramanujan constant $C_{LR} = 2^{-1/2}\prod_{p \equiv 3 \bmod 4}(1 - p^{-2})^{-1/2} = 0.76422365...$,
6) Twin primes constant $C_T = 2\prod_{p \geq 3}(1 - (p-1)^{-2}) = 1.3203236...$,
7) The special zeta values $\zeta(5)$, $\zeta(7)$, ...,
8) The imaginary parts of the critical zeros $\gamma_1 = 14.134725...$, $\gamma_2 = 21.022040...$, of the zeta function,

and many others.

A unified proof of the irrationality of the special values $L(n,\chi)$, $2 \leq n \in \mathbb{N}$, of beta *L*-function is put forward in this note. The first case of $n = 1$ seems to confirm that the Catalan constant $L(2,\chi)$ is an irrational number.

***Theorem* 1.** Let $\chi$ be the quadratic character, and let $L(s, \chi) = \sum_{n \geq 1} \chi(n) n^{-s}$ be a Dirichlet *L*-function. Then
1) The special values $L(n, \chi)^2$ are irrational numbers for all $2 \leq n \in \mathbb{N}$.
2) The special values $L(n, \chi)$ are irrational numbers for all $2 \leq n \in \mathbb{N}$.

Extensive information and wide range of views on the Catalan's constant $L(2, \chi)$ are given in [MW], [BL], [AV], [CR], [KW], [ST], [ZN], [WK], et al. The current research literature on the special values $L(2n, \chi)$ states the following:

(i) At least one of the seven numbers $L(2, \chi)$, $L(4, \chi)$, $L(6, \chi)$, $L(8, \chi)$, $L(10, \chi)$, $L(12, \chi)$, and $L(14, \chi)$ is an irrational number.

(ii) The value $L_2(2, \chi)$ of the 2-adic *L*-function $L_2(s, \chi)$ is an irrational number in $\mathbb{Q}_2 - \mathbb{Q}$.

The result (i) is derived from a hypergeometric analysis of the values $L(2n, \chi)$ and related constants, see [RL], and [ZN]; and result (ii) appears in [CR].

The second Section contains the basic material on the theory of irrationality and transcendence of real numbers; and the proof of Theorem 1 is given in the third Section.

## 2 Fundamental Concepts and Background

The basic concepts and results employed throughout this work are stated here as reference. A real number $\xi \in \mathbb{R}$ is called *rational* if $\xi = a/b$, where $a, b \in \mathbb{Z}$ are integers. Otherwise, the number is *irrational*. An irrational number $\xi \in \mathbb{R}$ is called *algebraic* if there exists a polynomial $f(x) \in \mathbb{Z}[x]$ such that $f(\xi) = 0$. Otherwise, the irrational number $\xi \in \mathbb{R}$ is called *transcendental*.

### 2.1 Continued Fractions

A continued fraction is a sequence of integers $[a_0, a_1, a_2, \ldots ]$, where $a_0 \in \mathbb{Z}$, and $a_i \in \mathbb{N}$ for $i \geq 1$. The set of continued fractions provides basefree representations of the real numbers. Rational numbers $a/b \in \mathbb{Q}$ have finite continued fractions $a/b = [a_0, a_1, a_2, \ldots, a_n]$, but irrational numbers $\xi \in \mathbb{R}$ have infinite continued fractions $\xi = [a_0, a_1, a_2, \ldots ]$.

The continued fractions of rational numbers are computed using the *Euclidean algorithm*, and the continued fractions of real algebraic numbers are computed using the *recursive quotient algorithm*, see [KM, p. 261]. Both of these algorithms use integer arithmetic to recursively generate the quotients $a_0, a_1, a_2, \ldots$ of any such real numbers. In contrast, except for some special cases, there is no known algorithm for computing the continued fractions of arbitrary transcendental numbers.

These algorithms are associated with the Gauss map

$$G(x) = \begin{cases} 0 & \text{if } x = 0, \\ 1/x - 1/[x] & \text{if } 0 < x \leq 1, \end{cases} \qquad (2)$$

where $[x]$ is the largest integer function. For a real number $\xi \in \mathbb{R}$, and an integer $k \geq 1$, the $k$th iterate $G^k(x) = [0, a_k, a_{k+1}, a_{k+2}, \ldots]$, of the Gauss map satisfies the following properties:

(1) If $\xi$ is a rational number, then $G^k(\xi) = 0$ for all large $k \geq k_0$.
(2) If $\xi$ is a quadratic irrational number, then $G^k(\xi) = \xi$ for all multiple $k = mt$ of the period $t \geq 1$.
(3) If $\xi$ is a nonquadratic irrational number, then $G^k(\xi) \neq \xi$ for all integers $k \geq 1$.

**Lemma 2.** (Fundamental recursion formulas) Let $\xi \in \mathbb{R}$ be a real number, let $[a_0, a_1, a_2, \ldots,]$ be its continued fraction. Then the following recursion relations hold.
(1) $p_{n+1} = a_{n+1} p_n + p_{n-1}, \quad p_0 = a_0, \quad p_1 = a_1 p_0 + 1,$
(2) $q_{n+1} = a_{n+1} q_n + q_{n-1}, \quad q_0 = 1, \quad q_1 = a_1,$
(3) $p_{n-1} q_n - p_n q_{n-1} = (-1)^n,$
(4) $p_{n-2} q_n - p_n q_{n-2} = (-1)^n a_n.$

As usual, the symbols $P(n)$ and $Q(n)$ denote the largest prime factor, and the largest squarefree factor respectively of the integer $n \geq 1$.

**Lemma 3.** For almost every algebraic irrational number $\xi \in \mathbb{R}$, the $n$th convergent $p_n/q_n$ has the following arithmetic properties.
(1) $P(q_n) \geq e^{\log q_n / \log \log q_n},$
(2) $Q(q_n) \geq q_n / (\log q_n)^{1-\delta}, \quad \delta > 0,$
(3) $\log q_n \to \pi^2 n / 12 \log 2$ as $n \to \infty$.

Confer [BU] for other arithmetic properties of the convergents, and other references.

**2.2 Criteria For Rationality And For Irrationality**

**Lemma 4.** (Criteria for rationality) If $\xi \in \mathbb{Q}$ is a rational number, then there exists a constant $c = c(\xi)$ such that

$$\left| \xi - \frac{p}{q} \right| \geq \frac{c}{q} \qquad (3)$$

holds for any rational fraction $p/q \neq \xi$.

This is a statement about the lack of effective or good approximations of an arbitrary rational number $\xi$ by other rational numbers. On the other hand, irrational numbers have effective approximations by rational numbers.

If the complementary inequality $\left| \xi - p/q \right| \leq c/q$ holds for infinitely many rational approximations $p/q$, then it already shows that the real number $\xi \in \mathbb{R}$ is almost irrational, so it is sufficient to prove the irrationality of real numbers. More precise results for testing the irrationality of an arbitrary real number are stated below.

***Lemma 5.*** (Criteria for irrationality) Let $\xi \in \mathbb{R}$ be a real number. If there exists a sequence of rational approximations $p_n/q_n$ such that $p_n/q_n \neq \xi$, and

$$\left| \xi - \frac{p_n}{q_n} \right| < \frac{1}{q_n^{1+\delta}} \tag{4}$$

for all $n \geq 1$ and some $\delta > 0$, then $\xi$ is an irrational number.

***Theorem 6.*** (Dirichlet) Let $\xi \in \mathbb{R}$ be a real number. If there exists a sequence of rational approximations $p_n/q_n$ such that $p_n/q_n \neq \xi$, and

$$\left| \xi - \frac{p_n}{q_n} \right| < \frac{1}{2q_n^2} \tag{5}$$

for all $n \geq 1$, then $\xi$ is an irrational number.

***Theorem 7.*** (Hurwitz) If $\xi \in \mathbb{R}$ is an irrational number, then there exists an infinite sequence of rational approximations $p_n/q_n$ such that

$$\left| \xi - \frac{p_n}{q_n} \right| < \frac{1}{\sqrt{5} q_n^2} \,. \tag{6}$$

Refer to [HW], [NV], et al, for more details.

The *p*-adic versions of these criteria are stated below; and an application to $L(2, \chi)$ and $\zeta(3)$ appears in [CR]. For a fixed prime $p \geq 2$, the symbol $\mathbb{Q}_p$ represents the set of *p*-adic real numbers, and the normalized *p*-adic absolute value on $\mathbb{Q}_p$ is defined by $\| x \|_p = 1/p$.

***Lemma 8.*** (Criteria of *p*-adic irrationality) Let $\xi \in \mathbb{Q}_p$ be a *p*-adic real number. If there exists some $\delta > 0$, and a sequence of rational approximations $p_n/q_n$ such that

$$\left\| \xi - \frac{p_n}{q_n} \right\|_p \leq \frac{1}{\left( \max \{ | p_n |, | q_n | \} \right)^{1+\delta}} \,, \tag{7}$$

where $p_n/q_n \neq \xi$, and $q_n \to \infty$ as $n \to \infty$, then $\xi$ is an irrational number.

***Theorem 9.*** (Ridout) Let $\xi \in \mathbb{R}$ be an irrational number. For every finite disjoint subsets of primes $S$ and $T$, and every real number $\varepsilon > 0$, the inequality

$$\left| \xi - \frac{p}{q} \right| \prod_{\ell \in S} | p |_\ell \prod_{\ell \in T} | p |_\ell < \frac{1}{q^{2+\varepsilon}} \tag{8}$$

has finitely many rational solutions $p/q$.

**Lemma 10.** Let $\xi \in \mathbb{R}$ be an irrational real number, and let $p/q < \xi < r/s$ be a pair of irreducible fractions such that $qr - ps = 1$. Then

$$\min\{ q^2(\xi - p/q),\ q^2(\xi - p/q) \} < 1/2. \tag{9}$$

The proof of this Lemma appears in [WT, p. 22].

**Theorem 11.** (Thue-Siegal-Roth)  Let $\varepsilon > 0$. An irrational number $\xi \in \mathbb{R}$ is algebraic if the inequality

$$\left| \xi - \frac{p}{q} \right| \leq \frac{1}{q^{2+\varepsilon}} \tag{10}$$

has only finitely many rational solutions $p/q \neq \xi$.

**Conjecture 12.** (Lange)  Let $\varepsilon > 0$. If $\xi \in \mathbb{R}$ is an algebraic irrational number, then the inequality

$$\left| \xi - \frac{p}{q} \right| \leq \frac{1}{q^2 (\log q)^{1+\varepsilon}} \tag{11}$$

has finitely many rational solutions $p/q \neq \xi$.

**Definition 13.** The irrationality measure $\mu(\xi)$ of a real number $\xi \in \mathbb{R}$ is defined by $\mu(\xi) = \limsup \{ \mu_0(\xi) : | \xi - p/q | < q^{-\mu_0(\xi)} \}$ for at most finitely many rational numbers $p/q$.

The definition of the irrationality measure is equivalent to the following: There is a constant $c > 0$ such that

$$\left| \xi - \frac{p}{q} \right| \geq \frac{c}{q^{\mu+\varepsilon}} \tag{12}$$

for all rational fractions $p/q \neq \xi$, and any $\varepsilon > 0$. For example, rational numbers have the irrationality measure $\mu(\xi) = 1$, and irrational numbers have the irrationality measure $\mu(\xi) \geq 2$.

**2.3 Autoconvolution of Dirichlet *L*-functions**
Some of the constants of unknown rationality or irrationality are defined by or have well known power series representations. In some cases these constants have Dirichlet series whose autoconvolutions are not very difficult to analyze. In the following result, the divisor function is defined by $d(n) = \#\{ d \mid n : 1 \leq d \leq n \}$.

**Lemma 14.** Let $\chi$ be the quadratic character, and let $L(s, \chi) = \sum_{n \geq 1} \chi(n) n^{-s}$ be a Dirichlet *L*-function. Then $L(s, \chi)^2 = \sum_{n \geq 1} d(n) \chi(n) n^{-s}$ for all complex numbers $s \in \mathbb{C}$ of definition $\mathfrak{Re}(s) > 1$.

*Proof*: The autoconvolution of $L(s, \chi)$ is

$$L(s,\chi)^2 = \sum_{k\geq 1}\frac{\chi(k)}{k^s}\sum_{n\geq 1}\frac{\chi(n)}{n^s} = \sum_{n\geq 1}\frac{a(n)}{n^s} = \sum_{n\geq 1}\frac{d(n)\chi(n)}{n^s}, \tag{13}$$

where $a(n) = \sum_{d|n}\chi(d)\chi(n/d) = \sum_{d|n}\chi(n) = d(n)\chi(n)$. This follows from the multiplicativity $\chi(a)\chi(b) = \chi(ab)$ of the character $\chi$, and the number of divisors $d(n) = \sum_{d|n} 1$ counting function. ∎

**Lemma 15.** Let $\varepsilon > 0$, and let $s \in \mathbb{C}$ be a complex number such that $a = \mathfrak{Re}(s) > 1$. Then, the partial sum $L(s,\chi)^2 = \sum_{n\geq 1} d(n)\chi(n) n^{-s}$ satisfies the followings estimates.
(1) $\sum_{n\leq x} d(n)\chi(n) n^{-s} = p_x / q_x$, where $p_x, q_x \in \mathbb{N}$, $c_0 x^{a-1-\varepsilon} \leq q_x \leq c_1 x^{a-1}$, and $c_0, c_1, c_2, c_3$ are constants, holds for almost every large integer $x > 0$.
(2) $\sum_{n\leq x} d(n)\chi(n) n^{-s} = p_x / q_x$, where $c_2 x^{a-\varepsilon} \leq q_x \leq c_3 x^a$, holds on a subset of integers $x \in D \subset \mathbb{N}$ of zero density in $\mathbb{N}$.

*Proof*: For simplicity take $s = 2$. Since $2 = \mathfrak{Re}(s) > 1$, the series is absolutely convergent, and the partial sum $\sum_{n\leq x} d(n)\chi(n) n^{-2}$ of $L(2,\chi)^2$ can be rearranged as

$$\sum_{n=1}^{x}\frac{\chi(n)d(n)}{n^2} = \sum_{n=0}^{x/4}\left(\frac{d(4n+1)}{(4n+1)^2} - \frac{d(4n+3)}{(4n+3)^2}\right)$$
$$= \sum_{n=0}^{x/4}\frac{(4n+3)^2 d(4n+1) - (4n+1)^2 d(4n+3)}{(4n+1)^2(4n+3)^2}. \tag{14}$$

This decomposition stems from the definition of the character $\chi(4n+1) = 1$, $\chi(4n+3) = -1$, and $\chi(2n) = 0$. Next, utilize order of the divisor function $d(n) = O(n^\varepsilon)$, $\varepsilon > 0$, and integral estimates, see Note 1, to arrive at

$$B_0 + \frac{c_0}{x} \leq \sum_{n=1}^{x}\frac{1}{n^2} \leq \sum_{n=1}^{x}\frac{\chi(n)d(n)}{n^2} \leq \sum_{n=1}^{x}\frac{n^\varepsilon}{n^2} \leq B_1 + \frac{c_1}{x^{1-\varepsilon}}, \tag{15}$$

where $B_0, B_1, c_0, c_1, c_2, \ldots, c_5$ are constants, for every sufficiently large number $x > 0$. The size of the set of integers for which (15) holds follows from the fact that the divisor function $d(4n+1) \neq d(4n+3)$ for almost every consecutive pair of integers $4n + 1$, $4n + 3$ as $n \to \infty$, see [PR]. Consequently, there are sufficiently many pairs such that $(4n + 3)^2 d(4n+1) - d(4n+3)(4n + 1)^2 > 0$ in the partial sum $\sum_{n\leq x} d(n)\chi(n) n^{-2} > 0$ for almost every sufficiently large $x > 0$.

As the partial sum is a finite sum of rational numbers, is it can be rewritten as a

$$\sum_{n=1}^{x}\frac{\chi(n)d(n)}{n^2} = \frac{p_x}{q_x}, \tag{16}$$

where $p_x, q_x \in \mathbb{N}$, $c_2 x^{a-1-\varepsilon} \leq q_x \leq c_3 x^{a-1}$, and $a = \mathfrak{Re}(s) = 2$, for all sufficiently large number $x > 0$. To confirm the second claim, note that

$$\sum_{n=1}^{x} \frac{\chi(n)d(n)}{n^2} = \sum_{n=0}^{x/4} \left( \frac{d(4n+1)}{(4n+1)^2} - \frac{d(4n+3)}{(4n+3)^2} \right) = \sum_{n=0}^{x/4} \frac{(16n+8)d(4n+1)}{(4n+1)^2 (4n+3)^2}, \qquad (17)$$

whenever $d(4n+1) = d(4n+3)$, and the order of the divisor function $d(n) = O(n^\varepsilon)$, $\varepsilon > 0$. Thus, the integral estimates yield

$$B_2 + \frac{c_4}{x^2} \leq n \sum_{k=0}^{x/4} \frac{(16n+8)d(4n+1)}{(4n+1)^2 (4n+3)^2} \leq B_3 + \frac{c_5}{x^{2-\varepsilon}}, \qquad (18)$$

where $B_2, B_3 > 0$ are constants, on a subset of integers $D = \#\{ n : d(4n+1) = d(4n+3) \}$ of zero density in $\mathbb{N}$. ∎

The distribution of repeated values $d(n) = d(n+k)$, $k \geq 1$, and other properties of the divisor function $d(n)$ are studied in [HB], [HD], and [PR].

**Note 1:** One of the possible integral approximations has the shape

$$\sum_{n=1}^{x} \frac{\chi(n)d(n)}{n^2} = \int_1^\infty \frac{dR(t)}{t^s}, \qquad (19)$$

where $R(x) = \sum_{n \leq x} \chi(n)d(n)$, this is estimated using standard analytic methods.

**Note 2:** The exact value of the partial sum is

$$\sum_{n=1}^{x} \frac{\chi(n)d(n)}{n^2} = \frac{N}{D} = \frac{N}{(1 \cdot 3 \cdot 5 \cdot 7 \cdot 9 \cdots x)^2}, \qquad (20)$$

where $N, D > 0$ are integers. However, the extremely large denominator $D$ is not suitable for proving the irrationality of the limiting sum, but the estimate (15) and (16) facilitate the determination of a suitable rational approximation for proving the irrationality of the limiting sum.

**2.4 The Irrationality of Some Constants**
The different analytical techniques utilized to confirm the irrationality, transcendence, and irrationality measures of many constants are important in the development of other irrationality proofs. Some of these results will be used later on.

***Theorem* 16.** The real numbers $\pi$, $\zeta(2) = \pi^2/6$, and $\zeta(3)$ are irrational numbers.

The proofs are widely available in the literature, refer to [HK], [BR], [AP], and others.

**Theorem 17.** For any fixed $n \in \mathbb{N}$, the followings statements are valid.

i) The real number $\zeta(2n) = (-1)^{n+1}(2\pi)^{2n} B_{2n} / 2(2n)!$ is an irrational number,

ii) The real number $L(2n+1, \chi) = (-1)^n \pi^{2n+1} E_{2n} / 2^{2n+2}(2n)!$ is an irrational number,

where $B_{2n}$ and $E_{2n}$ are the Bernoulli and Euler numbers respectively.

*Proof*: Apply the Lindemann-Weierstrass theorem to the transcendental number $\pi$. ∎

Theorems 1 and 17 completes the irrationality verifications of the special values $L(s, \chi)$ for any integer $s \geq 2$. Various other special values of the zeta function and *L*-functions are discussed in [CO] and similar sources.

## 3. Irrationality of The Special Values $L(2n, \chi)$ and $L(2n, \chi)^2$

Every irrational number $\xi \in \mathbb{R}$ has a representation as a *continued fraction power series*

$$\xi = a_0 + \sum_{k=1}^{\infty} \frac{(-1)^{n+1}}{q_n q_{n-1}}, \qquad (21)$$

where $q_n = a_n q_{n-1} + q_{n-2}$, $q_0 = 1$, $q_1 = a_1$ induced by the continued fraction $[a_0, a_1, a_2, \ldots]$ of $\xi$, see Lemma 2. In addition, this power series satisfies the stringent property $\left| \xi - a_0 - \sum_{k \leq n} (-1)^{-k} q_k^{-1} q_{k-1}^{-1} \right| = \left| \xi - p_n / q_n \right| < 1/2q_n^2$, see Theorems 7 and 8. The basic challenge of this approach to prove the irrationality of an arbitrary real number is the determination of the power series (21) or a closely related series that satisfies or nearly satisfies the stringent and ideal properties of the power series (21).

### 3.1 The Proof

Since the continued fractions of the special values $L(2n, \chi)$ and $L(2n, \chi)^2$ are unknown, the continued fraction power series (21) of these numbers are also unknown. However, pseudorandom nature of the *n*th term $d(n)\chi(n)/n^s$ of the special values of the Dirichlet series $L(s, \chi)^2$, see Lemma 14, makes it very closely related to the randomness of the *n*th term $(-1)^{n+1}/q_n q_{n-1}$ in the continued fraction (21) of the special values $L(2n, \chi)^2$.

**Theorem 1.** Let $\chi$ be the quadratic character, and let $L(s, \chi) = \sum_{n \geq 1} \chi(n) n^{-s}$ be a Dirichlet *L*-function. Then

1) The special values $L(n, \chi)^2$ are irrational numbers for all $2 \leq n \in \mathbb{N}$.
2) The special values $L(n, \chi)$ are irrational numbers for all $2 \leq n \in \mathbb{N}$.

*Proof*: Without loss in generality assume that $s = 2$, and consider $L(2, \chi)^2 = \sum_{n \geq 1} d(n)\chi(n) n^{-2}$, and its partial sum $\sum_{n \leq x} d(n)\chi(n) n^{-2}$. By Lemma 15, it has a rational approximation as

$$\sum_{n=1}^{x} \frac{\chi(n)d(n)}{n^2} = \frac{p_x}{q_x}, \tag{22}$$

where $p_x, q_x \in \mathbb{N}$, and $q_x \leq c_0 x$ for sufficiently large $x > 0$. Moreover, the alternating powers series $L(2, \chi)^2$ has the error term

$$0 < \left| L(2,\chi)^2 - \sum_{n=1}^{x} \frac{\chi(n)d(n)}{n^2} \right| = \left| \sum_{n>x}^{\infty} \frac{\chi(n)d(n)}{n^2} \right| \leq \frac{c_1}{x^{2-\varepsilon}}, \tag{23}$$

where $d(n) = O(n^\varepsilon)$, $\varepsilon > 0$, and $c_0, c_1, c_2, \ldots, c_5$ are constants. Now suppose that $L(2, \chi)^2 = A^2/B^2$ is a rational number. Then

$$0 < \frac{c_2}{x} \leq \frac{c_3}{q_x} \leq \left| L(2,\chi)^2 - \frac{p_x}{q_x} \right| \leq \frac{c_4}{x^{2-\varepsilon}} \leq \frac{c_5}{q_x^{2-\varepsilon}} \tag{24}$$

for infinitely many rational approximations $p_x/q_x \neq L(2,\chi)^2$ as $x \to \infty$, see Lemmas 4 and 5. Clearly, this is a contradiction, therefore $L(2, \chi)^2$ is an irrational number. ∎

### 3.2 Related Open Problems
Many interesting questions on the special values $L(2n, \chi)$ and $L(2n, \chi)^2$ of the beta *L*-function $L(s, \chi) = \sum_{n \geq 1} \chi(n) n^{-s}$ are open problems.

1) Are $L(2n, \chi)^2 = c\pi^d$, and $L(2n, \chi) = c\pi^d$, where $c, d \in \mathbb{Q}$ rational numbers?
2) Are $L(2n, \chi)^2$ and $L(2n, \chi)$ transcendental numbers?
3) What are the irrationality measures $\mu(L(2n, \chi))$ and $\mu(L(2n, \chi)^2)$?
4) Are $L(2n, \chi)^2$ and $L(2n, \chi)$ normal numbers?

Acknowledgement: Thanks to Professor Lima and his students for their comments.

Note on the Irrationality of *L*-Function Constants<:->